\documentclass[twoside]{amsart}
\usepackage{fullpage}
\usepackage{hyperref}
\begin{document}

\title{Liber Mathematicae: a web-based documentation and collaboration project 
for Mathematics} 
\author{Markus J. Pflaum}
\address{Department of Mathematics, University of Colorado, Boulder, CO 80309-0395, U.S.A.}
\email{Markus.Pflaum@colorado.edu}
\author{John Tuley}
\address{Department of Mathematics, University of Colorado, Boulder, CO 80309-0395, U.S.A.}
\email{John.Tuley@colorado.edu}

\thanks{presented at \textit{Workshop on the Future of Mathematics Journals}, MSRI,
        Berkeley, USA on February 15, 2011}

\maketitle
\thispagestyle{empty}
Traditionally, mathematical knowledge is published in printed media such as 
books or journals. With the advent of the Internet, a new method of publication 
became available. To date, however, most online mathematical publications do not 
employ the full capabilities of the medium. For example, the arXiv 
preprint server \cite{arXiv} hosts documents in Adobe's Portable Document Format (\textsc{PDF}), a 
format designed to transmit print documents between computers and maintain the 
formatting. Wikipedia \cite{wikipedia}, on the other hand, presents documents in-browser, 
but does not have support for presenting mathematics in scalable or copyable ways,
since all mathematics is presented as images, not text.

The languages of the World Wide Web, namely \textsc{HTML} and \textsc{XHTML}, are derived 
from Standardized Generalized Markup Language (\textsc{SGML}), originally, and Extensible 
Markup Language (\textsc{XML}), more recently. Unfortunately, \textsc{(X)HTML} has no 
provisions for displaying mathematics, a problem which the World Wide Web Consortium (W3C) 
solved with the \text{MathML} markup language. Like \textsc{XHTML}, \text{MathML} is a markup 
language with ``programs" stored as human-readable, plain-text files, formatted by a web 
browser or special viewer. As such, and as a standardized language, MathML is uniquely able to 
serve as a platform by which mathematics can be presented -- primarily on the Web -- but also 
processed by computer programs. In particular, computer algebra systems or tools for 
processing and disseminating text-based information such as search engines can 
handle mathematical content written in MathML. In an article in the \textit{Notices of the AMS}
from May 2005, \textsc{Miner} \cite{miner} explained  in detail such features and the importance
of MathML for the cummunication of mathematical content in the future. 

The real power of XML technology is that it allows to transcend a current 
``consumer culture" within the community of mathematicians, in which a few 
people produce content which is then consumed by many others, and reach what 
\textsc{Fischer} \cite{fischer} calls a 
``culture of participation" in which knowledge is created, shared, and acted upon by all 
parties. This goal can be achieved by encoding mathematics in browser-readable and 
easily-editable documents in MathML, instead of in read-only documents like \textsc{PDF}. 

The Liber Mathematicae project \cite{libermath} looks to bring the open source model of 
software development to mathematics publishing by employing cutting edge XML technology, 
high-quality mathematics  fonts for the Web from the STIX Fonts project \cite{stixfonts}
(see also \cite{newyorktimes}), and relational database 
technology to allow for a sophisticated version control and review process for the 
submitted mathematical content.  
We have developed a web site, where members of the mathematical community can not only view 
articles but can additionally participate in the creative process by contributing corrections,
suggest improvements, or by expanding on the original content. In contrast to traditional 
mathematics journals, the main goals of Liber Mathematicae are to have
articles which are expandable, correctable and dynamic, with tools for collaborative writing 
and open access to the entire mathematics community. 
Moreover due to their online nature, articles on Liber Mathematicae may contain more than 
static text and images and may in fact hold animations, live computational demonstrations, 
and so forth, and may use hyperlinking to strongly cross-reference other articles. An additional goal is to create a logical dependance tree for all 
mathematical theorems on Liber Mathematicae. We hope that with this new environment 
for communicating mathematical knowledge, the openness and cooperation will help to 
increase both the pace and quality of new mathematical research.

\bibliographystyle{alpha}

\end{document}